\documentclass[twocolumn,11pt,birkmult]{article}
\usepackage{times}

\newtheorem{definition}{Definition}
\newtheorem{remark}[definition]{Remark}

\newtheorem{theorem}[definition]{Theorem}

\newtheorem{corollary}[definition]{Corollary}

 \newtheorem{thm}{Theorem}[section]

 \newtheorem{defn}[thm]{Definition}
 
 \newtheorem{rem}[thm]{Remark}

\usepackage{amssymb}
\usepackage{amsmath}

\begin{document}
\global\def\refname{{\normalsize \it References:}}
\baselineskip 12.5pt
%
%
%

\title{Geometric Dynamics of Plasma in \\Jet Spaces with Berwald-Mo\'{o}r Metric\footnote{Proceedings of The 9th
International Conference on System Science and Simulation in Engineering (ICOSSSE-10),
Iwate, Japan, October 4-6, 2010}}

\date{}

\author{\hspace*{-10pt}
\begin{minipage}[t]{2.7in} \normalsize \baselineskip 12.5pt
\centerline{MIRCEA NEAGU}
\centerline{University Transilvania of Bra\c{s}ov}
\centerline{Faculty of Mathematics and Informatics}
\centerline{B-dul Iuliu Maniu, Nr. 50}
\centerline{BV 500091, Bra\c{s}ov, Romania}
\centerline{ mircea.neagu@unitbv.ro}
\centerline{mirceaneagu73@yahoo.com}
\end{minipage} \kern 0in
\begin{minipage}[t]{2.7in} \normalsize \baselineskip 12.5pt
\centerline{CONSTANTIN UDRI\c{S}TE}
\centerline{Univ. Politehnica of Bucharest}
\centerline{Faculty of Applied Sciences}
\centerline{Splaiul Independentei 313}
\centerline{060042 Bucharest, ROMANIA}
\centerline{ udriste@mathem.pub.ro}
\centerline{anet.udri@yahoo.com}
\end{minipage}
%
%
\\ \\ \hspace*{-10pt}
\begin{minipage}[b]{6.9in} \normalsize
\baselineskip 12.5pt {\it Abstract:}
In this paper we construct the differential equations of the stream lines
that characterize plasma regarded as a non-isotropic medium geometrized by
a jet rheonomic time-invariant Berwald-Mo\'{o}r metric.
Section 1 contains historical notes regarding the Plasma Physics and 
its geometrical description. 
Section 2 analyzes the generalized Lagrange geometrical approach of 
the non-isotropic
plasma on $1$-jet spaces. Section 3 studies the non-isotropic plasma as
 a medium geometrized by the jet
rheonomic Berwald-Mo\'{o}r metric.
\\ [4mm] {\it Key words and phrases:} jet Finsler spaces, relativistic rheonomic
Berwald-Mo\'{o}r metric, energy-stress-momentum d-tensor of non-isotropic
plasma, conservation laws, DEs of stream lines.
\\ [4mm] {\it Mathematics Subject Classification (2000):} 53B21, 53B40, 53C80.
\end{minipage}
\vspace{-10pt}}

\maketitle

\thispagestyle{empty} \pagestyle{empty}
%

%

\setcounter{equation}{0}
\section{Introduction}

\hspace{5mm}During that so-called the radiation epoch, in which photons are
strongly coupled with the matter, the interactions between the various
constituents of the Universal Matter include radiation-plasma coupling,
which is described by the plasma dynamics. Although it is not traditional to
characterize the radiation epoch by the dominance of plasma interactions,
however, it may be also called the plasma epoch (please, see \cite{Kleidis}).
This is because, in the plasma epoch, the electromagnetic interaction
dominates all the four fundamental physical forces (electrical, magnetic,
gravitational and nuclear).

In the present days, the Plasma Physics is an well established field of
Theoretical Physics, although the formulation of magnetohydrodynamics in a
curved space-time is a relatively new development (please see, Punsly \cite%
{Punsly}). The MHD processes in an isotropic space-time are intensively
studied by a lot of physicists. For example, the MHD equations in an
expanding Universe are investigated by Kleidis, Kuiroukidis, Papadopoulos
and Vlahos in \cite{Kleidis}. Considering the interaction of the
gravitational waves with the plasma in the presence of a weak magnetic
field, Papadopoulos also investigates the relativistic hydromagnetic
equations \cite{Papadopoulos}. The electromagnetic-gravitational dynamics
into plasmas with pressure and viscosity is studied by Das, DeBenedictis,
Kloster and Tariq in the paper \cite{Das}.

It is important to note that all preceding physical studies are done on an
isotropic four-dimensional space-time, represented by a semi- (pseudo-)
Riemannian space with the signature $(+,+,+,-)$. Consequently, the
Riemannian geometrical methods are used as a pattern over there.

From a geometrical point of view, using the Finlerian geometrical methods,
the plasma dynamics was extended on non-isotropic space-times by G\^{\i}r%
\c{t}u and Ciubotariu in the paper \cite{Girtu-Ciubotariu}. More
general, after the development of Lagrangian geometry of tangent bundle by
Miron and Anastasiei \cite{Mir-An}, the generalized Lagrange geometrical
objects describing the relativistic magnetized plasma were studied by M. G%
\^{\i}r\c{t}u, V. G\^{\i}r\c{t}u and Postolache in the paper \cite%
{Girtu-Girtu-Posto}.

According to Olver's opinion \cite{Olver}, we appreciate that the $1$-jet
spaces are basic geometrical objects in the study of classical and quantum
field theories. For such a reason, inspired by the geometrical methods of
geometric dynamics developed by Udri\c{s}te in \cite{Udr Geom Dyn}-\cite{UFO}, and
using as a pattern the Miron-Anastasiei's geometrical ideas from 
\cite{Mir-An}, Neagu recently developed in the book \cite{Neagu Carte} that
so-called "multitime Riemann-Lagrange geometry" on 1-jet
spaces, in the sense of distinguished d-connections, d-torsions, and
d-curvatures. We would like to point out that the geometrical construction
on $1$-jet spaces exposed in the monograph \cite{Neagu Carte} was initiated
by Asanov in \cite{Asanov[2]} and further developed by Neagu and Udri\c{s}%
te. Under the influence of the Riemann-Lagrange geometrical ideas from \cite%
{Neagu Carte}, the paper \cite{Neagu-Plasma MT} creates a multitime
extension on $1$-jet spaces of the geometrical objects that characterize
plasma in semi-Riemannian and Lagrangian approaches.

On the other hand, more studies of Russian researchers (Asanov \cite{Asanov[1]}%
, Mikhailov \cite{Mikhailov}, Garas'ko and Pavlov \cite{Garasko-Pavlov})
emphasize the importance in the study of physical fields of the Finsler
geometry which is characterized by the total equality of all non-isotropic
directions. For such a reason, Asanov, Pavlov and their co-workers underline
the important role played in the theory of space-time structure and
gravitation, as well as in unified gauge field theories, by the Berwald-Mo%
\'{o}r metric%
\begin{equation*}
F:TM\rightarrow \mathbb{R},\mathbb{\qquad }F(y)=\left(
y^{1}y^{2}...y^{n}\right) ^{\frac{1}{n}},
\end{equation*}%
whose Finsler geometry is intensively studied by Matsumoto and Shimada in
the paper \cite{Mats-Shimada}. Taking into account that our natural physical
intuition distinguishes four dimensions in a natural correspondence with the
material reality, in the framework of the $4$-dimensional linear space ($n=4$%
) with Berwald-Mo\'{o}r metric, Pavlov and his co-workers \cite%
{Garasko-Pavlov}, \cite{Pavlov} offer some new physical approaches and
geometrical interpretations for:

1. physical events = points in the 4-dimensional space;

2. straight lines = shortest curves;

3. intervals = distances between the points along of a straight line;

4. light pyramids $\Leftrightarrow $ light cones in a pseudo-Euclidean space.

In such a physical perspective and because of all preceding geometrical and
physical reasons, this paper is devoted to the development on the $1$-jet
space $J^{1}(\mathbb{R},M^{4})$ of the {\it geometric dynamics of plasma} endowed
with the \textit{relativistic rheonomic time-invariant Berwald-Mo\'{o}r metric}%
\begin{equation*}
\mathring{F}:J^{1}(\mathbb{R},M^{4})\rightarrow \mathbb{R},
\end{equation*}
\begin{equation*}
\mathring{F%
}(t,y)=\sqrt{h^{11}(t)}\,\,\sqrt[4]{y_{1}^{1}y_{1}^{2}y_{1}^{3}y_{1}^{4}},
\end{equation*}%
where $h_{11}(t)$ is a Riemannian metric on $\mathbb{R}$ and $%
(t,x^{1},x^{2},x^{3},x^{4},y_{1}^{1},y_{1}^{2},y_{1}^{3},y_{1}^{4})$ are the
coordinates on the $1$-jet space $J^{1}(\mathbb{R},M^{4})$. 

We underline that the geometry (in the sense of d-connections, d-torsions and
d-curvatures) produced by the jet rheonomic Berwald-Mo\'{o}r Finslerian
metric\textit{\ }$\mathring{F}$ is already completely done in the Neagu's
paper \cite{Neagu-B-M-4}. A part of this geometry
will be used in the next sections. Also, we mention that the sense of {\it geometric dynamics}
used in the present paper is different from those used by Udri\c ste \cite{Udr Geom Dyn}-\cite{UFO}.

\setcounter{equation}{0}
\section{Generalized Lagrange geometrical approach of the non-isotropic
plasma on $1$-jet spaces}

\hspace{5mm}Let $(\mathbb{R},h_{11}(t))$ be the set of real numbers endowed
with a Riemannian structure, where the coordinate $t$ plays the role of
relativistic time. The Christoffel symbol of the
Riemannian metric $h_{11}(t)$ is%
\begin{equation*}
\varkappa _{11}^{1}=\frac{h^{11}}{2}\frac{dh_{11}}{dt},\qquad h^{11}=\frac{1%
}{h_{11}}>0.
\end{equation*}

Let us consider that $M^{n}$ is a \textit{spatial real manifold} of
dimension $n$, whose local coordinates are $\left( x^{i}\right) _{i=%
\overline{1,n}}$. Notice that, in this Section, the latin indices run from $1
$ to $n$. Moreover, the Einstein convention of summation is used throughout
this work. Let $J^{1}(\mathbb{R},M)$ be the $1$-jet space of dimension $2n+1$%
, whose local coordinates are $(t,x^{i},y_{1}^{i})$. These transform by the
rules 
$$
\widetilde{t}=\widetilde{t}(t),\,\,\widetilde{x}^{i}=\widetilde{x}^{i}(x^{j}),\,\, 
\widetilde{y}_{1}^{i}={{\dfrac{\partial \widetilde{x}^{i}}{\partial x^{j}}}{%
\dfrac{dt}{d\widetilde{t}}}}\,\,\,y_{1}^{j}{,}%
$$
where $d\widetilde{t}/dt\neq 0$ and rank $(\partial \widetilde{x}%
^{i}/\partial x^{j})=n$.

Let $\mathcal{RGML}^{n}=\left( J^{1}(\mathbb{R},M),\text{ }%
G_{(i)(j)}^{(1)(1)}=h^{11}g_{ij}\right) $ be a relativistic rheonomic
generalized Lagrange space (for more details, please see Neagu \cite{Neagu
Carte}, \cite{Neagu-Rheon}), where $g_{ij}(t,x^{k},y_{1}^{k})$ is a metrical
d-tensor on $J^{1}(\mathbb{R},M)$, which is symmetrical, non-degenerate and
of constant signature.

Let us consider that $\mathcal{RGML}^{n}$ is endowed with a nonlinear
connection having the form \cite{Neagu Carte}, \cite{Neagu-Rheon}%
\begin{equation*}
\Gamma =\left( M_{(1)1}^{(i)}=-\varkappa _{11}^{1}y_{1}^{i},\text{ }%
N_{(1)j}^{(i)}\right) .
\end{equation*}%
The nonlinear connection $\Gamma $ produces on $J^{1}(T,M)$ the following
dual adapted bases of d-vectors and d-covectors:%
\begin{equation*}
\left\{ \frac{\delta }{\delta t},\frac{\delta }{\delta x^{i}},\frac{\partial 
}{\partial y_{1}^{i}}\right\} \subset \mathcal{X}(J^{1}(\mathbb{R}%
,M)),
\end{equation*}
\begin{equation*}
 \left\{ dt,dx^{i},\delta y_{1}^{i}\right\} \subset \mathcal{X}%
^{\ast }(J^{1}(\mathbb{R},M)),
\end{equation*}%
where%
\begin{equation*}
\frac{\delta }{\delta t}=\frac{\partial }{\partial t}+\varkappa
_{11}^{1}y_{1}^{m}\frac{\partial }{\partial y_{1}^{m}},\qquad \frac{\delta }{%
\delta x^{i}}=\frac{\partial }{\partial x^{i}}-N_{(1)i}^{(m)}\frac{\partial 
}{\partial y_{1}^{m}},
\end{equation*}%
\begin{equation*}
\delta y_{1}^{i}=dy_{1}^{i}-\varkappa
_{11}^{1}y_{1}^{i}dt+N_{(1)m}^{(i)}dx^{m}.
\end{equation*}%
It is important to note that the d-tensors on the $1$-jet space $J^{1}(%
\mathbb{R},M)$ behave like classical tensors. For example, on the $1$-jet
space $J^{1}(\mathbb{R},M)$ we have the global metrical d-tensor%
\begin{equation*}
\mathbb{G}=h_{11}dt\otimes dt+g_{ij}dx^{i}\otimes dx^{j}+h^{11}g_{ij}\delta
y_{1}^{i}\otimes \delta y_{1}^{j},
\end{equation*}%
which is endowed with the physical meaning of non-isotropic gravitational
potential. Obviously, the d-tensor $\mathbb{G}$ has the adapted components%
\begin{equation*}
\mathbb{G}_{AB}=\left\{ 
\begin{array}{llll}
h_{11}, & \text{for} & A=1, & B=1\medskip  \\ 
g_{ij}, & \text{for} & A=i, & B=j\medskip  \\ 
h^{11}g_{ij}, & \text{for} & A=_{(i)}^{(1)}, & B=_{(j)}^{(1)}\medskip  \\ 
0, & \text{otherwise.} &  & 
\end{array}%
\right. 
\end{equation*}

Following the geometrical ideas of Asanov \cite{Asanov[2]} and Neagu \cite%
{Neagu-Rheon}, the preceding geometrical ingredients lead us to the the
Cartan canonical $\Gamma $-linear connection (given in adapted components)%
\begin{equation*}
C\Gamma =\left( \varkappa _{11}^{1},\text{ }G_{j1}^{k},\text{ }L_{jk}^{i},%
\text{ }C_{j(k)}^{i(1)}\right) ,
\end{equation*}%
where%
\begin{equation}
\begin{array}{c}
\medskip G_{j1}^{k}=\dfrac{g^{km}}{2}\dfrac{\delta g_{mj}}{\delta t},\\
L_{jk}^{i}=\dfrac{g^{im}}{2}\left( \dfrac{\delta g_{jm}}{\delta x^{k}}+%
\dfrac{\delta g_{km}}{\delta x^{j}}-\dfrac{\delta g_{jk}}{\delta x^{m}}%
\right) , \\ 
C_{j(k)}^{i(1)}=\dfrac{g^{im}}{2}\left( \dfrac{\partial g_{jm}}{\partial
y_{1}^{k}}+\dfrac{\partial g_{km}}{\partial y_{1}^{j}}-\dfrac{\partial g_{jk}%
}{\partial y_{1}^{m}}\right) .%
\end{array}
\label{Cartan-Multi-Time}
\end{equation}%
In the sequel, the Cartan linear connection $C\Gamma $, given by (\ref%
{Cartan-Multi-Time}), induces the $\mathbb{R}$-horizontal ($h_{\mathbb{R}}-$%
) covariant derivative%
\begin{eqnarray*}
D_{1k(1)(l).../1}^{1i(j)(1)...} &\!\!\!\!=\!\!\!\!&\frac{\delta D_{1k(1)(l)...}^{1i(j)(1)...}%
}{\delta t}+D_{1k(1)(l)...}^{1i(j)(1)...}\varkappa
_{11}^{1}\\
&\!\!\!\!+\!\!\!\!&D_{1k(1)(l)...}^{1m(j)(1)...}G_{m1}^{i}+D_{1k(1)(l)...}^{1i(m)(1)...}G_{m1}^{j}\\
&\!\!\!\!+\!\!\!\!&D_{1k(1)(l)...}^{1i(j)(1)...}\varkappa _{11}^{1}+... \\
&\!\!\!\!-\!\!\!\!&D_{1k(1)(l)...}^{1i(j)(1)...}\varkappa_{11}^{1}-D_{1m(1)(l)...}^{1i(j)(1)...}G_{k1}^{m} \\
&\!\!\!\!-\!\!\!\!&D_{1k(1)(l)...}^{1i(j)(1)...}\varkappa
_{11}^{1}-D_{1k(1)(m)...}^{1i(j)(1)...}G_{l1}^{m}...,
\end{eqnarray*}%
the $M$-horizontal ($h_{M}-$) covariant derivative%
\begin{eqnarray*}
D_{1k(1)(l)...|p}^{1i(j)(1)...}&\!\!\!\!=\!\!\!\!&\frac{\delta D_{1k(1)(l)...}^{1i(j)(1)...}%
}{\delta x^{p}}%
+D_{1k(1)(l)...}^{1m(j)(1)...}L_{mp}^{i}\\
&\!\!\!\!+\!\!\!\!&D_{1k(1)(l)...}^{1i(m)(1)...}L_{mp}^{j}+...\\
&\!\!\!\!-\!\!\!\!&D_{1m(1)(l)...}^{1i(j)(1)...}L_{kp}^{m}-D_{1k(1)(m)...}^{1i(j)(1)...}L_{lp}^{m}-...
\end{eqnarray*}%
and the vertical ($v-$) covariant derivative%
\begin{eqnarray*}
D_{1k(1)(l)...}^{1i(j)(1)...}|_{(p)}^{(1)} &\!\!\!\!=\!\!\!\!&\frac{\partial
D_{1k(1)(l)...}^{1i(j)(1)...}}{\partial y_{1}^{p}}%
+D_{1k(1)(l)...}^{1m(j)(1)...}C_{m(p)}^{i(1)}\\
&\!\!\!\!+\!\!\!\!&D_{1k(1)(l)...}^{1i(m)(1)...}C_{m(p)}^{j(1)}+...\\
&\!\!\!\!-\!\!\!\!&D_{1m(1)(l)...}^{1i(j)(1)...}C_{k(p)}^{m(1)}\\
&\!\!\!\!-\!\!\!\!&D_{1k(1)(m)...}^{1i(j)(1)...}C_{l(p)}^{m(1)}-...,
\end{eqnarray*}%
where%
\begin{eqnarray*}
D&\!\!\!\!=\!\!\!\!&D_{1k(1)(l)...}^{1i(j)(1)...}(t,x^{r},y_{1}^{r})\frac{\delta }{\delta t}%
\otimes \frac{\delta }{\delta x^{i}}\otimes \frac{\partial }{\partial
y_{1}^{j}}\otimes dt \\
&\!\!\!\!\otimes\!\!\!\!& dx^{k}\otimes \delta y_{1}^{l}\otimes ...
\end{eqnarray*}%
is an arbitrary d-tensor on $J^{1}(\mathbb{R},M)$.

\begin{remark}
The Cartan covariant derivatives produced by $C\Gamma $ have the
metrical properties%
\begin{equation*}
\begin{array}{lll}
h_{11/1}=h_{\text{ \ }/1}^{11}=0, & h_{11|k}=h_{\text{ \ }|k}^{11}=0,\medskip \\ 
h_{11}|_{(k)}^{(1)}=h^{11}|_{(k)}^{(1)}=0,\medskip \\ 
g_{ij/1}=g_{\text{ \ }/1}^{ij}=0, & g_{ij|k}=g_{\text{ \ }|k}^{ij}=0, \medskip \\ 
g_{ij}|_{(k)}^{(1)}=g^{ij}|_{(k)}^{(1)}=0.%
\end{array}%
\end{equation*}
\end{remark}

For the study of the magnetized non-viscous plasma dynamics, in a
relativistic generalized Lagrangian geometrical approach on 1-jet spaces, we
use the following geometrical objects \cite{Girtu-Girtu-Posto}, \cite%
{Neagu-Plasma MT}:

\begin{enumerate}
\item the unit relativistic time dependent velocity-d-field of a test
particle, which is given by%
\begin{equation*}
U=u_{1}^{i}(t,x^{k},y_{1}^{k})\frac{\partial }{\partial y_{1}^{i}},
\end{equation*}%
where, if we take $\varepsilon ^{2}=h^{11}g_{ij}y_{1}^{i}y_{1}^{j}>0$, then
we put $u_{1}^{i}=y_{1}^{i}/\varepsilon $. Obviously, we have $%
h^{11}u_{i1}u_{1}^{i}=1$, where $u_{i1}=g_{im}u_{1}^{m}$;

\item the distinguished relativistic time dependent $2$-form of the
(electric field)-(magnetic induction), which is given by%
\begin{equation*}
H=H_{ij}(t,x^{k},y_{1}^{k})dx^{i}\wedge dx^{j};
\end{equation*}

\item the distinguished relativistic time dependent $2$-form of the
(electric in\-duc\-tion)-(magnetic field), which is given by%
\begin{equation*}
G=G_{ij}(t,x^{k},y_{1}^{k})dx^{i}\wedge dx^{j};
\end{equation*}

\item the relativistic time dependent Minkowski energy d-tensor of the
electromagnetic field inside the non-isotropic plasma, which is given by%
\begin{eqnarray*}
E&\!\!\!\!=\!\!\!\!&E_{ij}(t,x^{k},y_{1}^{k})dx^{i}\otimes
dx^{j}\\
\medskip \\ 
&\!\!\!\!+\!\!\!\!&h^{11}E_{ij}(t,x^{k},y_{1}^{k})\delta y_{1}^{i}\otimes \delta
y_{1}^{j}.
\end{eqnarray*}%
The adapted components of the relativistic time dependent Minkowski energy
are defined by%
\begin{equation*}
E_{ij}=\dfrac{1}{4}g_{ij}H_{rs}G^{rs}+g^{rs}H_{ir}G_{js},
\end{equation*}%
where $G^{rs}=g^{rp}g^{sq}G_{pq}$. Moreover, we suppose that the adapted
components of the relativistic time dependent Min\-kow\-ski energy verify
the \textit{non-isotropic Lorentz con\-di\-tions} \cite{Neagu-Plasma MT}%
\begin{equation}
E_{i|m}^{m}u_{1}^{i}=0,\qquad E_{i}^{m}|_{(m)}^{(1)}u_{1}^{i}=0,
\label{Lorentz-conditions-Multi-Time}
\end{equation}%
where $E_{i}^{m}=g^{mp}E_{pi}$. Obviously, if we use the notations $%
H_{r}^{m}=g^{mp}H_{pr}$ and $G_{i}^{r}=g^{rs}G_{si}$, we have%
\begin{equation*}
E_{i}^{m}=\frac{1}{4}\delta _{i}^{m}H_{rs}G^{rs}-H_{r}^{m}G_{i}^{r},
\end{equation*}%
where $\delta _{i}^{m}$ is the Kronecker symbol.
\end{enumerate}

In our jet generalized Lagrangian geometrical approach, the relativistic
time dependent non-isotropic plasma is characterized by an
energy-stress-momentum d-tensor $\mathcal{T}$, which is defined by \cite%
{Girtu-Girtu-Posto}, \cite{Neagu-Plasma MT}%
\begin{eqnarray*}
\mathcal{T}&\!\!\!\!=\!\!\!\!&\mathcal{T}_{ij}(t,x^{k},y_{1}^{k})dx^{i}\otimes dx^{j}\\
&\!\!\!\!+\!\!\!\!&h^{11}\mathcal{T}_{ij}(t,x^{k},y_{1}^{k})\delta y_{1}^{i}\otimes \delta y_{1}^{j},
\end{eqnarray*}%
where \cite{Das}%
\begin{equation}
\mathcal{T}_{ij}=\left( \mathbf{\rho }+\dfrac{\mathbf{p}}{c^{2}}\right)
h^{11}u_{i1}u_{j1}+\mathbf{p}g_{ij}+E_{ij}.  \label{stress-Multi-Time-1}
\end{equation}%
The entities $c=$ constant, $\mathbf{p}=\mathbf{p}(t,x^{k},y_{1}^{k})$ and $%
\mathbf{\rho }=\mathbf{\rho }(t,x^{k},y_{1}^{k})$ have the physical meanings
of: - the speed of light, the non-isotropic hydrostatic pressure and the
non-isotropic proper mass density of plasma. Notice that the adapted
components of the energy-stress-momentum d-tensor $\mathcal{T}$, that
characterizes the non-isotropic plasma, are given by%
\begin{equation}
\mathcal{T}_{CF}=\left\{ 
\begin{array}{llll}
\mathcal{T}_{ij}, & \text{for} & C=i, & F=j\medskip  \\ 
h^{11}\mathcal{T}_{ij}, & \text{for} & C=_{(i)}^{(1)}, & F=_{(j)}^{(1)}%
\medskip  \\ 
0, & \text{otherwise.} &  & 
\end{array}%
\right.   \label{stress-Multi-Time-2}
\end{equation}

In the jet generalized Lagrange framework for plasma, we postulate that the
following \textit{non-isotropic conservation laws} of the components (\ref%
{stress-Multi-Time-1}) and (\ref{stress-Multi-Time-2}) are true:%
\begin{equation}
\mathcal{T}_{A:M}^{M}=0,\qquad \forall \text{ }A\in \left\{ 1,\text{ }i,%
\text{ }_{(i)}^{(1)}\right\} ,  \label{cons-general-M-T}
\end{equation}%
where the capital latin letters $A,M,...$ are indices of kind $1,$ $i$ or $%
_{(i)}^{(1)}$, "$_{:M}$" represents one of the local covariant derivatives $%
h_{\mathbb{R}}-$, $h_{M}-$ or $v-$ and%
\begin{equation*}
\mathcal{T}_{A}^{M}=\mathbb{G}^{MD}\mathcal{T}_{DA}
=\left\{ 
\begin{array}{ll}
\mathcal{T}_{i}^{m}, & \text{for } A=i, \,\, M=m  \\ 
\mathcal{T}_{i}^{m}, & \text{for } A=_{(i)}^{(1)}, \,\, M=_{(m)}^{(1)}
\\ 
0, & \text{otherwise.}   
\end{array}%
\right. 
\end{equation*}%
Note that the d-tensor $\mathcal{T}_{i}^{m}$ is given by the formula%
\begin{equation*}
\mathcal{T}_{i}^{m}=g^{mp}\mathcal{T}_{pi}=\left( \mathbf{\rho }+\dfrac{%
\mathbf{p}}{c^{2}}\right) h^{11}u_{1}^{m}u_{i1}+\mathbf{p}\delta
_{i}^{m}+E_{i}^{m}.
\end{equation*}

It is easy to see that the jet non-isotropic conservation laws (\ref%
{cons-general-M-T}) reduce to the following local non-isotropic conservation
equations:%
\begin{equation}
\mathcal{T}_{i|m}^{m}=0,\qquad \mathcal{T}_{i}^{m}|_{(m)}^{(1)}=0.
\label{cons-eq-M-T}
\end{equation}%
Moreover, by direct computations, we deduce that the non-isotropic
conservation equations (\ref{cons-eq-M-T}) become%
\begin{equation}
\begin{array}{l}
h^{11}\left[ \left( \mathbf{\rho }+\dfrac{\mathbf{p}}{c^{2}}\right) u_{1}^{m}%
\right] _{|m}u_{i1}+\left( \mathbf{\rho }+\dfrac{\mathbf{p}}{c^{2}}\right)
h^{11}u_{1}^{m}u_{i1|m}\medskip  \\ 
+\mathbf{p}_{,,i}-g_{ir}\overset{h}{\mathcal{F}^{r}}%
=0,\medskip  \\ 
h^{11}\left. \left[ \left( \mathbf{\rho }+\dfrac{\mathbf{p}}{c^{2}}\right)
u_{1}^{m}\right] \right\vert _{(m)}^{(1)}u_{i1}+\left( \mathbf{\rho }+\dfrac{%
\mathbf{p}}{c^{2}}\right) h^{11}u_{1}^{m}\medskip  \\ 
\cdot u_{i1}|_{(m)}^{(1)}+\mathbf{p}_{\#(i)}^{\text{ \ }(1)}-g_{ir}\overset{v}{\mathcal{F}}\text{ }%
\!\!^{r1}=0,%
\end{array}
\label{cons-local-M-T}
\end{equation}%
where $\mathbf{p}_{,,i}=\delta \mathbf{p}/\delta x^{i}$, $\mathbf{p}%
_{\#(i)}^{\text{ \ }(1)}=\partial \mathbf{p}/\partial y_{1}^{i}$ and

\begin{itemize}
\item $\overset{h}{\mathcal{F}^{r}}=-g^{rs}E_{s|m}^{m}$ is the \textit{%
non-isotropic horizontal Lorentz force};

\item $\overset{v}{\mathcal{F}}$ $\!\!^{r1}=-g^{rs}E_{s}^{m}|_{(m)}^{(1)}$
is the \textit{non-isotropic vertical} \textit{Lorentz d-tensor force.}
\end{itemize}

Contracting now the non-isotropic conservation equations (\ref%
{cons-local-M-T}) with $u_{1}^{i}$ and taking into account the non-isotropic
Lorentz conditions (\ref{Lorentz-conditions-Multi-Time}), we find the 
\textit{non-isotropic continuity equations} of plasma, namely%
\begin{equation}
\begin{array}{l}
\left[ \left( \mathbf{\rho }+\dfrac{\mathbf{p}}{c^{2}}\right) u_{1}^{m}%
\right] _{|m}+\mathbf{p}_{,,m}u_{1}^{m}=0,\medskip  \\ 
\left. \left[ \left( \mathbf{\rho }+\dfrac{\mathbf{p}}{c^{2}}\right)
u_{1}^{m}\right] \right\vert _{(m)}^{(1)}+\mathbf{p}_{\#(m)}^{\text{ \ \ }%
(1)}u_{1}^{m}=0,%
\end{array}
\label{Continuity-Lagrange}
\end{equation}%
where we also used the equalities%
\begin{equation*}
\begin{array}{l}
0=h^{11}u_{i1}u_{1|m}^{i}=\dfrac{1}{2}\left(
h^{11}u_{i1}u_{1}^{i}\right) _{,,m}=-h^{11}u_{i1|m}u_{1}^{i}, \\ 
\medskip0=h^{11}u_{i1}u_{1}^{i}|_{(m)}^{(1)}=\dfrac{1}{2}\left(
h^{11}u_{i1}u_{1}^{i}\right) _{\#(m)}^{\text{ \ \ }%
(1)}\\
=-h^{11}u_{i1}|_{(m)}^{(1)}u_{1}^{i},%
\end{array}
\end{equation*}
the symbols "$_{,,m}$" and "$_{\#(m)}^{\text{ \ \ }(1)}$" being the
derivative operators $\delta /\delta x^{m}$ and $\partial /\partial y_{1}^{m}$.

Replacing the continuity laws (\ref{Continuity-Lagrange}) into the
conservation equations (\ref{cons-local-M-T}), we find the \textit{%
non-isotropic relativistic Euler equations} for plasma, namely%
\begin{equation}
\begin{array}{l}
\medskip \left( \mathbf{\rho }+\dfrac{\mathbf{p}}{c^{2}}\right)
h^{11}u_{i1|m}u_{1}^{m}-\mathbf{p}_{,,m}\left( h^{11}u_{1}^{m}u_{i1}-\delta
_{i}^{m}\right) \\
-g_{im}\overset{h}{\mathcal{F}^{m}}=0, \\ 
\left( \mathbf{\rho }+\dfrac{\mathbf{p}}{c^{2}}\right)
h^{11}u_{i1}|_{(m)}^{(1)}u_{1}^{m}-\mathbf{p}_{\#(m)}^{\text{ \ \ }%
(1)}\left( h^{11}u_{1}^{m}u_{i1}-\delta _{i}^{m}\right) \\
-g_{im}\overset{v}{%
\mathcal{F}}\!\text{ }\!^{m1}=0.%
\end{array}
\label{Euler-Lagrange}
\end{equation}

If we take now $y_{1}^{m}=dx^{m}/dt$, then we have%
\begin{equation*}
u_{1}^{m}=\frac{1}{\varepsilon _{0}}\frac{dx^{m}}{dt}=\frac{dx^{m}}{ds}, 
\end{equation*}
\begin{equation*}
\varepsilon _{0}^{2}=h^{11}(t)g_{ij}(t,x,dx/dt)\frac{dx^{i}}{dt}%
\frac{dx^{j}}{dt},
\end{equation*}%
where $s$ is a natural parameter of the curve $c=(x^{k}(t))$, having the
property $ds/dt=\varepsilon _{0}$. Introducing this $u_{1}^{m}$ into the
non-isotropic Euler equations (\ref{Euler-Lagrange}), we obtain the
equations of the \textit{non-isotropic stream lines} for jet plasma, which
are given by the following second order DE systems:

\begin{itemize}
\item \textit{horizontal} non-isotropic stream line DEs:%
\begin{eqnarray*}
\dfrac{d^{2}x^{k}}{ds^{2}}&\!\!\!\!+\!\!\!\!&\left[ L_{rm}^{k}-\dfrac{c^{2}}{\mathbf{p}+%
\mathbf{\rho }c^{2}}\delta _{r}^{k}\mathbf{p}_{,,m}\right] \dfrac{dx^{r}}{ds}%
\dfrac{dx^{m}}{ds}\\
&\!\!\!\!=\!\!\!\!&\dfrac{h_{11}c^{2}}{\mathbf{p}+\mathbf{\rho }c^{2}}\left[ 
\overset{h}{\mathcal{F}^{k}}\!\!-\!\!g^{km}\mathbf{p}_{,,m}\right] \!\!+\!\!\dfrac{N_{(1)m}^{(k)}}{\varepsilon _{0}}\dfrac{dx^{m}}{ds}\\
&\!\!\!\!-\!\!\!\!&\dfrac{%
h^{11}N_{(1)m}^{(p)}g_{pr}}{\varepsilon _{0}}\dfrac{dx^{r}}{ds}\dfrac{dx^{m}%
}{ds}\dfrac{dx^{k}}{ds} \\
&\!\!\!\!-\!\!\!\!&\dfrac{h^{11}N_{(1)m}^{(r)}}{2}\dfrac{\partial g_{pq}}{\partial y_{1}^{r}}%
\dfrac{dx^{p}}{ds}\dfrac{dx^{q}}{ds}\dfrac{dx^{m}}{ds}\dfrac{dx^{k}}{ds};%
\end{eqnarray*}

\item \textit{vertical} non-isotropic stream line DEs:%
$$\left[C_{r(m)}^{k(1)}-\dfrac{c^{2}}{\mathbf{p}+\mathbf{\rho }c^{2}}\delta
_{r}^{k}\mathbf{p}_{\#(m)}^{\text{ \ \ }(1)}\right] \dfrac{dx^{r}}{ds}\dfrac{%
dx^{m}}{ds}$$
$$\begin{array}{l}
=\dfrac{h_{11}c^{2}}{\mathbf{p}+\mathbf{\rho }c^{2}}\left[ 
\overset{v}{\mathcal{F}}\!\text{ }\!^{k1}-g^{km}\mathbf{p}_{\#(m)}^{\text{ \
\ }(1)}\right]\medskip  \\ 
+\dfrac{h^{11}}{2}\dfrac{\partial g_{pq}}{\partial y_{1}^{r}}\dfrac{dx^{p}}{%
ds}\dfrac{dx^{q}}{ds}\dfrac{dx^{r}}{ds}\dfrac{dx^{k}}{ds}.%
\end{array}$$
\end{itemize}

\begin{remark}
If the metrical d-tensor $g_{ij}(t,x,y)$ is Finslerian-like one, that is
we have%
\begin{equation*}
g_{ij}(t,x,y)=\frac{h_{11}}{2}\frac{\partial ^{2}F^{2}}{\partial
y_{1}^{i}\partial y_{1}^{j}},
\end{equation*}%
where $F:J^{1}(\mathbb{R},M)\rightarrow \mathbb{R}_{+}$ is a jet Finslerian
metric, then we use the \textbf{canonical spatial nonlinear connection} $%
N=\left( N_{(1)j}^{(k)}\right) $ of the jet Finsler space, whose general
formula is given in \cite{Neagu-Rheon}. Consequently, the DEs of the \textbf{%
stream lines} of plasma in non-isotropic jet Finsler spaces reduce to:

\begin{itemize}
\item \textbf{horizontal} non-isotropic stream line DEs:%
\begin{equation*}
\begin{array}{l}
\dfrac{d^{2}x^{k}}{ds^{2}}+\left[ L_{rm}^{k}-\dfrac{c^{2}}{\mathbf{p}+%
\mathbf{\rho }c^{2}}\delta _{r}^{k}\mathbf{p}_{,,m}\right] \dfrac{dx^{r}}{ds}%
\dfrac{dx^{m}}{ds}\\
=\dfrac{h_{11}c^{2}}{\mathbf{p}+\mathbf{\rho }c^{2}}\left[ 
\overset{h}{\mathcal{F}^{k}}-g^{km}\mathbf{p}_{,,m}\right]+\dfrac{N_{(1)m}^{(k)}}{\varepsilon _{0}}\dfrac{dx^{m}}{ds} \\ 
\medskip-\dfrac{%
h^{11}N_{(1)m}^{(p)}g_{pr}}{\varepsilon _{0}}\dfrac{dx^{r}}{ds}\dfrac{dx^{m}%
}{ds}\dfrac{dx^{k}}{ds};%
\end{array}%
\end{equation*}

\item \textbf{vertical} non-isotropic stream line DEs:%
\begin{equation}
\mathbf{p}_{\#(m)}^{\text{ \ \ }(1)}\left[ h_{11}g^{mk}-\dfrac{dx^{m}}{ds}%
\dfrac{dx^{k}}{ds}\right] =h_{11}\overset{v}{\mathcal{F}}\text{ \ }%
\!\!\!^{k1},  \label{vertical-Minkowski}
\end{equation}%
where $\varepsilon _{0}=F$ and$,$ if the generalized Christoffel symbols of $%
g_{ij}(t,x,y)$ are%
\begin{equation*}
\Gamma _{jk}^{i}(t,x,y)=\frac{g^{im}}{2}\left( \frac{\partial g_{jm}}{%
\partial x^{k}}+\frac{\partial g_{km}}{\partial x^{j}}-\frac{\partial g_{jk}%
}{\partial x^{m}}\right) ,
\end{equation*}%
then we have%
$$\begin{array}{l}
N_{(1)l}^{(k)}y_{1}^{l}=\Gamma _{pq}^{k}y_{1}^{p}y_{1}^{q}+\frac{\varkappa
_{11}^{1}g^{km}}{2}\frac{\partial g_{pq}}{\partial x^{m}}y_{1}^{p}y_{1}^{q}\medskip\\
\medskip+\frac{g^{km}}{2}\frac{\partial g_{mp}}{\partial t}y_{1}^{p}-\frac{\varkappa
_{11}^{1}}{2}y_{1}^{k}.
\end{array}$$

\hspace{5mm}More particular, if we have a jet Minkowski-like metric $F=F(y)$%
, then the \textbf{horizontal} non-isotropic stream line DEs of plasma
simplify as:%
\begin{equation}
\begin{array}{l}
\dfrac{d^{2}x^{k}}{ds^{2}}+\left[ L_{rm}^{k}-\dfrac{c^{2}}{\mathbf{p}+%
\mathbf{\rho }c^{2}}\delta _{r}^{k}\mathbf{p}_{,,m}\right] \dfrac{dx^{r}}{ds}%
\dfrac{dx^{m}}{ds}=\medskip \\ 
=\dfrac{h_{11}c^{2}}{\mathbf{p}+\mathbf{\rho }c^{2}}\left[ \overset{h}{%
\mathcal{F}^{k}}-g^{km}\mathbf{p}_{,,m}\right] .%
\end{array}
\label{horizontal-Minkowski}
\end{equation}
\end{itemize}
\end{remark}


\setcounter{equation}{0}
\section{The non-isotropic plasma as a medium geometrized by the jet
rheonomic Berwald-Mo\'{o}r metric}

\hspace{5mm}Let us note that in this Section we work with a fixed 4-dimensional spatial
manifold $M^{4}$. Consequently, throughout this Section, the
latin indices run only from $1$ to $4$. Also note that in this Section we
will focus only on the \textit{jet relativistic rheonomic Berwald-Mo\'{o}r
metric} \cite{Neagu-B-M-4}%
\begin{equation}
\mathring{F}(t,y)=\sqrt{h^{11}(t)}\,\, \sqrt[4]{%
y_{1}^{1}y_{1}^{2}y_{1}^{3}y_{1}^{4}}.  \label{rheon-B-M}
\end{equation}

Using the notation $G_{1111}=y_{1}^{1}y_{1}^{2}y_{1}^{3}y_{1}^{4},$ the
fundamental metrical d-tensor produced by the relativistic rheonomic
Berwald-Mo\'{o}r metric (\ref{rheon-B-M}) is given by (no sum by $i$ or $j$)

$$g_{ij}(t,x,y)=\frac{h_{11}(t)}{2}\frac{\partial ^{2}\mathring{F}^{2}}{%
\partial y_{1}^{i}\partial y_{1}^{j}}$$
\begin{equation}
=\frac{\left( 1-2\delta _{ij}\right) 
\sqrt{G_{1111}}}{8}\frac{1}{y_{1}^{i}y_{1}^{j}},  \label{g-jos-(ij)}
\end{equation}
where $\delta _{ij}$ is the Kronecker symbol. The matrix $g=(g_{ij})$ admits
the inverse $g^{-1}=(g^{jk})$, whose entries are (no sum by $j$ or $k$)%
\begin{equation}
g^{jk}=\frac{2(1-2\delta ^{jk})}{\sqrt{G_{1111}}}y_{1}^{j}y_{1}^{k}\text{.}
\label{g-sus-(jk)}
\end{equation}

Using some general formulas from the paper \cite{Neagu-Rheon}, we find the
following geometrical results \cite{Neagu-B-M-4}:

\begin{theorem}
For the jet relativistic rheonomic Berwald-Mo\'{o}r metric (\ref{rheon-B-M}%
), the \textit{energy action functional}%
\begin{equation*}
\mathbb{\mathring{E}}(x(\cdot))=\int_{a}^{b}\sqrt{%
y_{1}^{1}y_{1}^{2}y_{1}^{3}y_{1}^{4}}\,\, h^{11}\sqrt{h_{11}}\,\,dt,
\end{equation*}%
where $y_{1}^{i}=dx^{i}/dt$, produces on the $1$-jet space $J^{1}(\mathbb{R}%
,M^{4})$ the \textbf{canonical time dependent spray} 
\begin{equation*}
\mathcal{\mathring{S}}=\left( H_{(1)1}^{(i)}=-\frac{\varkappa _{11}^{1}}{2}%
y_{1}^{i},\text{ }G_{(1)1}^{(i)}=-\frac{\varkappa _{11}^{1}}{3}%
y_{1}^{i}\right) .
\end{equation*}
\end{theorem}

\begin{corollary}
The \textbf{canonical nonlinear connection} on the jet space of first order $%
J^{1}(\mathbb{R},M^{4}),$ associated to the jet relativistic rheonomic
Berwald-Mo\'{o}r me\-tric (\ref{rheon-B-M}), is given by%
\begin{equation}
\mathring{\Gamma}=\left\{ \begin{array}{l} M_{(1)1}^{(i)}=2H_{(1)1}^{(i)}=-\varkappa
_{11}^{1}y_{1}^{i},\\
\text{ }N_{(1)j}^{(i)}=\frac{G_{(1)1}^{(i)}}{\partial
y_{1}^{j}}=-\frac{\varkappa _{11}^{1}}{3}\delta _{j}^{i}.
\end{array}
\right.
\label{nlc-B-M}
\end{equation}

\end{corollary}

\begin{remark}
The nonlinear connection (\ref{nlc-B-M}) produces the dual \textit{adapted
bases} of d-vector fields%
\begin{equation}
\left\{ \frac{\delta }{\delta t},\text{ }\frac{\delta 
}{\delta x^{i}},\text{ }\dfrac{\partial }{\partial
y_{1}^{i}}\right\} \subset \mathcal{X}(E)  \label{a-b-v}
\end{equation}%
where,
$$
\begin{array}{l}
 \frac{\delta }{\delta t}=\frac{\partial }{\partial t}+\varkappa
_{11}^{1}y_{1}^{p}\frac{\partial }{\partial y_{1}^{p}},\medskip\\
\medskip\text{ }\frac{\delta 
}{\delta x^{i}}=\frac{\partial }{\partial x^{i}}+\frac{\varkappa _{11}^{1}}{3%
}\frac{\partial }{\partial y_{1}^{i}}.\\
\end{array}$$
and d-covector fields%
\begin{equation}
\left\{ dt,\text{ }dx^{i},\delta y_{1}^{i}\right\} \subset 
\mathcal{X}^{\ast }(E),  \label{a-b-co}
\end{equation}%
where
$$\begin{array}{l}
\text{ }\delta y_{1}^{i}=dy_{1}^{i}-\varkappa
_{11}^{1}y_{1}^{i}dt-\frac{\varkappa _{11}^{1}}{3}dx^{i}
\end{array}$$
and $E=J^{1}(\mathbb{R},M^{4})$. Note that the distinguished geometrical
elements of the adapted bases (\ref{a-b-v}) and (\ref{a-b-co}) transform
like classical tensors.
\end{remark}

On the $1$-jet space $J^{1}(\mathbb{R},M^{4})$, we will describe the Cartan
canonical connection produced by the relativistic rheonomic Berwald-Mo\'{o}r
metric (\ref{rheon-B-M}) in local adapted components. Thus, using the
formulas (\ref{Cartan-Multi-Time}), (\ref{g-jos-(ij)}) and (\ref{g-sus-(jk)}%
), by direct computations, we obtain \cite{Neagu-B-M-4}:

\begin{theorem}
The Cartan canonical $\mathring{\Gamma}$-linear connection$,$ produced by
the rheonomic Berwald-Mo\'{o}r metric (\ref{rheon-B-M})$,$ has the following
adapted local components:%
\begin{equation}
C\mathring{\Gamma}=\left( \varkappa _{11}^{1},\text{ }G_{j1}^{k}=0,\text{ }%
L_{jk}^{i}=\frac{\varkappa _{11}^{1}}{3}C_{j(k)}^{i(1)},\text{ }%
C_{j(k)}^{i(1)}\right) ,  \label{Cartan-Berwald-Moor}
\end{equation}%
where, if we use the notation%
\begin{equation*}
A_{jk}^{i}=\frac{2\delta _{j}^{i}+2\delta _{k}^{i}+2\delta _{jk}-8\delta
_{j}^{i}\delta _{jk}-1}{8}
\end{equation*}%
(no sum by $i$, $j$ or $k$), then we have%
\begin{equation*}
C_{j(k)}^{i(1)}=A_{jk}^{i}\cdot \frac{y_{1}^{i}}{y_{1}^{j}y_{1}^{k}}\text{
(no sum by }i,\text{ }j\text{ or }k\text{).}
\end{equation*}
\end{theorem}

\begin{remark}
The below properties of the d-tensor $C_{j(k)}^{i(1)}$ are true (sum by $m$):%
\begin{equation*}
C_{j(k)}^{i(1)}=C_{k(j)}^{i(1)},\quad C_{j(m)}^{i(1)}y_{1}^{m}=0,\quad
C_{j(m)}^{m(1)}=0.
\end{equation*}%
For more details, please see also the papers \cite{At-Bal-Neagu} and \cite%
{Mats-Shimada}.
\end{remark}

\begin{remark}
The coefficients $A_{ij}^{l}$ have the following values:%
\begin{equation*}
A_{ij}^{l}=\left\{ 
\begin{array}{ll}
-\dfrac{1}{8}, & i\neq j\neq l\neq i\medskip  \\ 
\dfrac{1}{8}, & i=j\neq l\text{ or }i=l\neq j\text{ or }j=l\neq i\medskip 
\\ 
-\dfrac{3}{8}, & i=j=l.%
\end{array}%
\right. 
\end{equation*}
\end{remark}

The geometric dynamics of the non-isotropic plasma regarded as a medium
geometrized by the jet rheonomic Berwald-Mo\'{o}r metric (\ref{rheon-B-M}%
) is obtained using the DEs of stream lines (\ref{horizontal-Minkowski}) and
(\ref{vertical-Minkowski}) for the particular Berwald-Mo\'{o}r geometrical
objects (\ref{g-jos-(ij)}), (\ref{g-sus-(jk)}), (\ref{a-b-v}) and (\ref%
{Cartan-Berwald-Moor}). Consequently, we find the following geometrical
DEs for non-isotropic plasma:

\begin{itemize}
\item the \textit{horizontal} Berwald-Mo\'{o}r non-isotropic stream line DEs
are given by%
\begin{eqnarray*}\dfrac{d^{2}x^{k}}{ds^{2}}&+&\dfrac{c^{2}}{\mathbf{p}+\mathbf{\rho }c^{2}}%
\mathbf{p}_{,,m}\left( 1-4\delta ^{km}\right) \dfrac{dx^{m}}{ds}\dfrac{dx^{k}%
}{ds}\\
&=&\dfrac{h_{11}c^{2}}{\mathbf{p}+\mathbf{\rho }c^{2}}\overset{h}{%
\mathcal{F}^{k}},
\end{eqnarray*}
where $k\in \left\{ 1,2,3,4\right\} $ is a fixed index and we do sum by $m$;

\item the \textit{vertical} Berwald-Mo\'{o}r non-isotropic stream line DEs
are given by%
\begin{equation*}
\mathbf{p}_{\#(m)}^{\text{ \ \ }(1)}\left( 1-4\delta ^{mk}\right) \dfrac{%
dx^{m}}{ds}\dfrac{dx^{k}}{ds}=h_{11}\overset{v}{\mathcal{F}}\text{ \ }%
\!\!\!^{k1},
\end{equation*}%
where $k\in \left\{ 1,2,3,4\right\} $ is a fixed index and we do sum by $m$.
\end{itemize}

\begin{remark}
In the particular case when the hydrostatic pressure is dependent only by $t$
and $x$ (i.e. we have $\mathbf{p}=\mathbf{p}(t,x^{k})$)$,$ the DEs of stream
lines for non-isotropic plasma endowed with Berwald-Mo\'{o}r metric become:

\begin{itemize}
\item \textbf{horizontal} Berwald-Mo\'{o}r non-isotropic stream line DEs:%
\begin{eqnarray*}
\dfrac{d^{2}x^{k}}{ds^{2}}&+&\dfrac{c^{2}}{\mathbf{p}+\mathbf{\rho }c^{2}}%
\mathbf{p}_{,m}\left( 1-4\delta ^{km}\right) \dfrac{dx^{m}}{ds}\dfrac{dx^{k}%
}{ds}\\
&=&\dfrac{h_{11}c^{2}}{\mathbf{p}+\mathbf{\rho }c^{2}}\overset{h}{%
\mathcal{F}^{k}},
\end{eqnarray*}%
where $\mathbf{p}_{,m}=\partial \mathbf{p}/\partial x^{m};$

\item \textbf{vertical} Berwald-Mo\'{o}r non-isotropic stream line DEs:%
\begin{equation*}
\overset{v}{\mathcal{F}}\text{ \ }\!\!\!^{k1}=0.
\end{equation*}
\end{itemize}
\end{remark}

\textbf{Open Problem. }There exist real physical interpretations for the
preceding jet Fin\-sler-Berwald-Mo\'{o}r geometric dynamics of non-isotropic
plasma?


\end{document}